\documentclass[12pt,reqno]{amsart}

\usepackage{graphicx}

\newtheorem{prop}{Proposition}
\newtheorem{defin}{Definition}
\newtheorem{thm}{Theorem}

\newtheorem{corol}{Corollary}
\newtheorem{lem}{Lemma}

\newtheorem{rem}{Remark}

\newcommand{\h}{{\mathfrak h}}

\newcommand{\N}{{\mathbb N}}

\newcommand{\R}{{\mathbb R}}\newcommand{\RR}{{\mathbb R}^2}
\newcommand{\pr}{\R\cup\infty }  

\newcommand{\x}{{\mathcal X}}
\newcommand{\Z}{{\mathbb Z}}

\newcommand{\bo}{\partial} 
\newcommand{\bom}{\bo\Om} 
\newcommand{\const}{\mbox{const}} 
\newcommand{\con}{\mbox{const}} 
\newcommand{\crit}{\mbox{crit}} 

\newcommand{\el}{\text{Ell}}  
\newcommand{\fix}{\text{fix}}  

\newcommand{\hy}{\text{Hyp}}  
\newcommand{\id}{\text{Id}}   

\newcommand{\mac}{\text{SL}(2,\R)}   
\newcommand{\psl}{\text{PSL}(2,\R)}   
\newcommand{\pgl}{\text{PGL}(2,\R)}   
\newcommand{\gl}{\text{GL}(2,\R)}   
\newcommand{\neut}{\text{Dgt}}  
\newcommand{\per}{\text{Per}}  
\newcommand{\pa}{\text{Par}}  

\newcommand{\sph}{S^2} 
\newcommand{\vol}{\text{vol}}  

\newcommand{\tO}{\tilde{O}}

\newcommand{\tpsi}{\tilde{\psi}}
\newcommand{\tPsi}{\tilde{\Psi}}

\newcommand{\al}{\alpha}
\newcommand{\be}{\beta}
\newcommand{\ga}{\gamma}\newcommand{\Ga}{\Gamma}
\newcommand{\ep}{\varepsilon}
\newcommand{\la}{\lambda}
\newcommand{\ka}{\kappa}

\newcommand{\Om}{\Omega}

\begin{document}

\bibliographystyle{plain}

\title[Periodic billiard orbits]
{A few remarks on periodic orbits for planar billiard tables}

\author{Eugene Gutkin}
\address{
IMPA, Estrada Dona Castorina 110, Rio de Janeiro, Brasil 22460-320
}
\email{gutkin@impa.br}


\date{\today}

\begin{abstract} 
I announce a solution of the conjecture about the  measure of periodic points 
for planar billiard tables.
The theorem says that if $\Om\subset\R^2$ is
a compact domain with piecewise $C^3$ boundary, then the set of periodic
orbits for the billiard in $\Om$ has measure zero.
Here I outline a proof. A complete version will appear elsewhere. 
\end{abstract}

\maketitle

\section{Introduction}             \label{intro}
Let $\Om\subset\RR$ be a connected, compact domain with a piecewise differentiable
boundary. We turn $\Om$ into a billiard table by letting a dimensionless ball
move freely inside, bouncing elastically off of $\bo\Om$. This dynamical system
is the {\em billiard in $\Om$}. It is a very special case of the geodesic flow.
The flow lines are the billiard orbits in  $\Om$; their properties reflect the
geometry of the domain. The simplicity of the setting makes the planar billiard
an ideal ``playground": We can ask meaningful questions about
the billiard dynamics essentially from scratch.\footnote{This is a free quote of 
the opinion of G.D. Birkhoff \cite{Bi}.
}
The billiard in a smooth and
convex domain  is the ``billiard ball problem" of G.D. Birkhoff \cite{Bi};
his study of periodic orbits for the billiard ball problem remains 
a crowning achievement in geometric dynamics.

The properties of billiard orbits continue to attract our attention,
suggesting a variety of questions, which seem very basic and turn
out to be out of reach \cite{Gut03}. Many of the open questions concern periodic 
billiard orbits. These orbits correspond to closed, inscribed polygons $P\subset\Om$ 
that have the following property:
At every corner of $P$ the two angles between $P$ and $\Om$ are equal.
Thus, periodic billard orbits correspond to {\em harmonic polygons} \cite{Bi}.
Using variational considerations, Birkhoff discovered that every smooth, convex billiard 
table had at least two infinite families of  distinct harmonic polygons, yielding the famous Birkhoff
periodic orbits.

Most open problems in this subject ask if billiard 
tables have ``sufficiently many" periodic orbits.\footnote{
For instance, it is not known if every convex polygon has a periodic orbit
\cite{Gut03}. The question is open even for obtuse triangles.
}
Going in the opposite direction, it is natural to ask how big the set of  periodic 
billiard orbits may be.  
There are examples of billiard tables where periodic orbits 
admit continuous deformations. A classical example is the round disc: Rotating harmonic
polygons, we obtain one-parameter families of periodic orbits of the same kind. 

The property of having continuous families of periodic orbits holds 
for an arbitrary ellipse;\footnote{
The disc being a special case of the ellipse.
} 
there are other examples of billiard tables with this
property. However, all of these examples yield only one-parameter periodic
families. The question arises if there exists a billiard table with 
a two-parameter periodic family. In such a table we could shoot the billiard ball, 
varying independently
the initial position and the initial direction, and it
would keep returning to the same position and direction.
The intuition says that no billiard table has the magic property.
In fact, a well known conjecture states that no billiard table has a two-parameter family of
periodic orbits \cite{Gut01,Gut03}.

\vspace{4mm}

Periodic billiard orbits are important not
only in dynamics; they are important in analysis.
With a compact euclidean domain $\Om\subset\R^d$ we associate the infinite sequence 
$0\le\la_1\le\la_2\le\cdots$ of eigenvalues of a laplacean in $\Om$. 
There are two classical eigenvalue problems, corresponding to the dirichlet and the neumann
laplaceans. The dirichlet eigenfunctions vanish on  $\bo\Om$, while the values of neumann
eigenfunctions are arbitrary on  $\bo\Om$, but their normal derivatives
vanish. In what follows we will refer to $\la_i$ as the neumann or dirichlet eigenvalues 
of $\Om$.
The motivation for studying the asymptotics of eigenvalues
came from  physics, and the mathematical formulation goes back to D. Hilbert.
Seminal results about 
the eigenvalues of euclidean domains were obtained by H. Weyl
around 1911-1915 \cite{WeA}.

Let $N_{\Om}(\la)$ be the number of eigenvalues of $\Om\subset\R^d$ (either neumann or dirichlet)
that are less than or equal to $\la^2$. Weyl's law is the asymptotic formula
\begin{equation}   \label{weyl_eq}
N_{\Om}(\la) = c_1\vol(\Om)\la^d + o(\la^d).
\end{equation}
The constant $c_1$ is known and depends only on the dimension.
Weyl's formula has numerous extensions. It is valid, mutatis mutandis,
on closed riemannian manifolds, riemannian manifolds with boundaries,
for higher order differential operators, etc \cite{SaVa}.

The formula  equation~\eqref{weyl_eq} holds under mild assumptions, 
but it yields only the leading term of the asymptotics of
eigenvalues. Note that the leading term does not depend
at all on the boundary $\bo\Om$.  Weyl conjectured that
the asymptotic expansion equation~\eqref{weyl_eq} 
continues, and that the second term involves $\vol(\bo\Om)$.
The precise formulation of his conjecture is the equation  
\begin{equation}   \label{ivri_eq}
N_{\Om}(\la) = c_1\vol(\Om)\la^d \pm c_2\vol(\bo\Om)\la^{d-1} + o(\la^{d-1}).
\end{equation}
The constant $c_2$ is also known and is determined by $d$; 
the sign of the second term depends on the boundary condition. 

In \cite{Iv80} V. Ivrii proved Weyl's conjecture for compact, $C^{\infty}$  
domains $\Om\subset\R^d$, provided that billiard orbits in $\Om$ satisfy an assumption.
In order to prove equation~\eqref{ivri_eq},  Ivrii had to assume that periodic
orbits for the billiard in $\Om$ spanned a set of measure zero.\footnote{For instance, with respect to
the canonical invariant measure.}

During 1979/80 Ivrii had repeatedly asked the Moscow
billiard community whether the assumption was redundant.
The experts assured him that no billiard table had a set of positive measure of
periodic orbits; however, they could not supply a proof.\footnote{
See http://www.math.toronto.edu/ivrii/Research/preprints/GradTalk3.pdf
for a colorful account of Ivrii's interactions with participants 
of Sinai's seminar.} 

Ever since the paper \cite{Iv80} appeared, the billiard community at large 
felt that every euclidean billiard table
satisfied Ivrii's assumption. Lacking a proof, it became a conjecture.
In what follows I will refer to it as the {\em measure of periodic points conjecture},
or simply the {\em periodic points conjecture}.
The conjecture and its variants figure prominently in the dynamics literature 
\cite{Tab,Gut01,Gut03,R89}, and
in the analysis literature \cite{SaVa,PeSt}.
A few partial results
towards the conjecture have been obtained, and I will comment on them below. 

\medskip

In this work I announce a positive solution of the measure of periodic points conjecture 
for planar  billiard tables.
More precisely, the theorem says that if $\Om\subset\R^2$ is
a compact, piecewise $C^3$ domain, then the measure of the set of periodic
orbits for the billiard in $\Om$ is zero.
This paper is preliminary, and I only outline the  proofs.
A complete version will appear elsewhere. 

\medskip

From here on, $\Om$ is a planar domain, unless  
explicitly stated otherwise. I will comment now on the
variants of the conjecture that are in the literature. 
The billiard dynamics in $\Om$ has two realisations: As the billiard flow,
and as the billiard map.
The phase space of the billiard flow is three-dimensional.
Phase points are the pairs (position, velocity) where the position can be anywhere in $\Om$.
The billiard map is the poincare map with respect to
the cross-section that consists of phase points where the position is in $\bo\Om$.
This cross-section is the phase space of the billiard map.
The liouville measure is invariant for the flow, and induces
the canonical measure for the map. The flow version and the map version
of the periodic points conjecture are equivalent.
In what follows we will work only with the billiard map. 

The phase space for the billiard map is a compact surface with a boundary\footnote{
If $\Om$ is simply connected, then it is a topological annulus.
}
that consists of fixed points; we will ignore the boundary,
and denote by $\Psi=\Psi(\Om)$  the set of interior points.
The set $\per(\Om)$ of periodic points is a countable union of the sets of
$n$-periodic points, $\per_n(\Om)$, where $2\le n<\infty$. A seemingly different variant of the
conjecture says that for every $n$ the set $\per_n(\Om)\subset\Psi$
is nowhere dense. This claim is formally weaker than the equation $\mu(\per_n(\Om))=0$;
it turns out, however, that the two statements are equivalent. See sections~\ref{main},~\ref{crit}.

\medskip

There are three kinds of partial results towards the periodic points conjecture. 
First, there are results that establish the claim under additional
assumptions on the billiard table; e. g., assuming the (piecewise) real
analyticity of $\bo\Om$ \cite{SaVa}. 
It is worth pointing out that in this case the claim is obvious: 
If the derivative of a real analytic function vanishes
on an open set, then the function is constant. Thus, if a real analytic billiard table
$\Om$ has an open set of $n$-periodic orbits, then the $n$th iterate of the billiard
map for $\Om$ is the identity map, which is absurd.
Results of the second kind prove the claim for generic (smooth) billiard tables 
\cite{PeSt}. It is obvious that having continuous families of periodic orbits
is a non-generic property. In fact, for the generic domain $\Om\subset\R^d$ 
the sets $\per_n(\Om)$ are finite \cite{PeSt}.
Hence, generically, $\mu(\per(\Om))=0$. 

Since $\per(\Om)=\cup_{n\in\N}\per_n(\Om)$, establishing the 
periodic points conjecture amounts to
proving the claim for every $n$. Partial results of the third kind do this for small $n$.
Since $2$-periodic points are contained in the set of phase points orthogonal to $\bom\subset\R^d$,
the smallest period for which the claim is not obvious is $3$. 
In \cite{R89} M. Rychlik proved that the set of $3$-periodic points for a planar, convex,
infinitely differentiable billiard table has empty interior.
His proof is based on the following proposition.\footnote{
I will refer to it as Rychlik's lemma in what follows.
}  
For $s\in\bo\Om$ let 
$\Psi_s\subset\Psi$ be the set of phase points corresponding to shooting the billiard ball
from point $s$ in arbitrary directions. Rychlik's lemma states that $\per_3(\Om)\cap\Psi_s$
is nowhere dense in $\Psi_s$, for any $s\in\bo\Om$. Rychlik's approach is based on 
symbolic calculations involving
the length of orbits, and he relied on a software package to perform them.
In \cite{St} L. Stojanov pointed out that the lengths of periodic orbits
do not change under deformations; this remark disposed of the computer part in
Rychlik's argument. In \cite{Vo97} Ya.B. Vorobets proved Rychlik's lemma
using elementary geometry. His proof requires only the differentiability of
$\bo\Om$; moreover, it works for euclidean domains in any dimension. In 
\cite{W94} M.P. Wojtkowski gave another proof of Rychlik's lemma; 
it is based on an analysis of Jacobi fields for the billiard.
In \cite{BZh} Y. Baryshnikov and V. Zharnitsky prove that $\per_3(\Om)$ is nowhere dense
using external differentiable systems (EDS).\footnote{According to a private communication of 
Zharnitsky, they are working on $\per_4(\Om)$ using the EDS approach.
}

\section{Overview of the proof}        \label{over}
To simplify the exposition, we restrict the discussion in the technical part of the paper
to convex billiard tables. We remove the restriction in section~\ref{main}, where we
state our main results, and outline the proofs.
In the discussion that follows we assume that $\Om$ is convex. 

Let $\psi$ be a $n$-periodic phase point. Suppose that $\psi$ is contained in a ball
$O$ of $n$-periodic points. Then the differential $DF^n$ is identically one on $O$.
In order to study the equation $DF_{\psi}^n=1$, we pass to the projectivised differential,
$PDF$, that acts on the  projectivised tangent bundle of $\Psi$. See section~\ref{mirror}.

Elements of the projectivised tangent bundle have a geometric optics interpretation, and
the action of $PDF$ is given by the {\em mirror equation}. See equation~\eqref{mirr_eq}.
An orbit $\psi_0,\psi_1=F(\psi_0),\ldots,\psi_{n-1}=F^{n-1}(\psi_0),\ldots$
yields a sequence of fractional-linear transformations on $\R\cup\infty$.
A number $x_0\in\R\cup\infty$ is interpreted as the signed 
focusing  distance of an infitesimal
light beam sent along the orbit. Then $x_i=PDF^i\circ x_0$ is the focusing distance for
the beam that we obtain after $i$ reflections in $\bo\Om$. The transformation 
$x_n=PDF^n\circ x_0$ determines the type of our periodic point. We have $DF_{\psi}^n=\pm 1$,
iff $x_n\equiv x_0$. In this case we say that $\psi$ is a degenerate periodic point.

Fractional-linear transformations $x_i=x_i(x_0)$ involve products of matrices. 
These expressions are awkward to work with. We get around this difficulty by working instead 
with the derivatives $dx_i/dx_0$.  
See equation~\eqref{der_foc_eq}. Analysis of these derivatives 
yields consequences for the geometry of harmonic polygons corresponding to
degenerate periodic points. See Proposition~\ref{neutr_geom_prop}.

In particular, we obtain an identity relating the distance between consecutive
reflection points of the billiard ball and the geometry of $\bo\Om$ at these points. 
See equation~\eqref{basic_geom_eq}. 

Suppose now that we have an open ball $O\subset\Psi$ of $n$-periodic points.
Let $s_0,s_1$ be the endpoints of the directed chord in $\Om$ corresponding 
to a phase point $\psi\in O$. Note that each point $s_i$ varies
independently in an open arc, $\Ga_i\subset\bo\Om$.
Equation~\eqref{basic_geom_eq} yields a functional identity, equation~\eqref{imposs_eq},
relating the radii of curvature of $\bo\Om$ at $s_0,s_1$.

In section~\ref{calcul} we study this functional identity. 
Under the assumption of continuous differentiability of
the radius of curvature for $\bo\Om$, we prove that equation~\eqref{imposs_eq}
has no solutions. See Proposition~\ref{gener_anal_prop} and Corollary~\ref{no_iden_cor}.
   
In section~\ref{main} we obtain the main results of the paper.
In this section  we extend the discussion 
of sections~\ref{mirror},~\ref{per},~\ref{calcul} to compact planar
domains with a piecewise $C^3$ boundary, disposing of the convexity assumption.
 
Theorem~\ref{emt_inter_thm} states that 
the set of periodic points for the billiard in $\Om$ has an empty interior.
Periodic points in an open set are degenerate; the converse
may fail, e. g., the billiard in $\Om$ may have isolated degenerate periodic points.
Theorem~\ref{odd_neut_type_thm} says that degenerate periodic points cannot have odd periods.
Theorem~\ref{main_conv_thm} is about the measure of the set of
periodic points. It says that $\mu(\per(\Om))=0$.
The proof combines Theorem~\ref{emt_inter_thm} with Theorem~\ref{per_crit_thm}
which is the main result of section~\ref{crit}.

\medskip

Thus, Theorem~\ref{main_conv_thm} removes the measure zero assumption
in Ivrii's theorem \cite{Iv80}. By Theorem~\ref{main_conv_thm} and Ivrii's result, 
Weyl's conjecture holds for all compact planar domains with $C^{\infty}$ boundary. 

\medskip

In section~\ref{crit} we show that $\per(\Om)$
is nowhere dense iff $\mu(\per(\Om))=0$. We point out that 
this is a special feature of the billiard dynamics, by exhibiting
infinitely differentiable, conservative diffeomorphisms of the $2$-torus that do not
satisfy this property.

Periodic billiard orbits are given by a variational
principle. 
This fact and Theorem~\ref{emt_inter_thm}
imply  Theorem~\ref{per_crit_thm}
about the structure of the set of periodic billiard orbits. It says that 
every $\per_n(\Om)\subset\Psi$  is a finite union of curves and points.

The billiard map for a compact, convex, planar domain, $\Om$, is an area preserving twist map.  
It is known that the claim $\mu(\per(\Om))=0$ does not extend to such maps \cite{R89}.
Hence, it is not surpising that our approach is based on
the geometry of planar billiards, and that we dispose of the convexity
assumption. Does the claim $\mu(\per(\Om))=0$ extend to general nonplanar,
two-dimensional billiard tables? Section~\ref{rem} addresses this question. 

Let $\Om$ be a compact riemannian surface with a (piecewise) smooth boundary.
The billiard map for $\Om$ is well defined, and let $\mu$ be the canonical invariant measure
on its phase space $\Psi=\Psi(\Om)$.
Let $\per(\Om)\subset\Psi$ be the set of periodic points for the billiard in $\Om$.
Is $\mu(\per(\Om))=0$? Suppose that $\bom$ is geodesically convex. Is $\mu(\per(\Om))=0$ then?
In section~\ref{rem} we give negative answers to these questions.
Let $\sph$  be the round sphere. We construct infinitely smooth,
convex billiard tables $\Om\subset\sph$ such that the set $\per(\Om)$ has
a nonempty interior. Moreover, we can make the domains $\Om$ in these examples  
arbitrarily small. See Corollary~\ref{sph_ex_cor}. We can also dispose of the convexity assumption.
Thus, $\mu(\per(\Om))=0$  is not a universal billiard property.

\section{The projectivised tangent bundle}        \label{set}
We begin by establishing the basic notation. Recall that the billiard table is a compact 
domain $\Om\subset\RR$ with a piecewise $C^3$ boundary $\bo\Om$.
Each connected component of $\bo\Om$ is positively oriented, and parametrised by
arclength.
Let $\Psi=\bo\Om\times[0,\pi]$ be the phase space of the billiard map, $F:\Psi\to\Psi$.
In the standard coordinatization,
phase points are the pairs $\psi=(s,\al)\in\Psi$, where $s\in\bom$ is the
footpoint, and $0\le\al\le\pi$ is the direction angle.
The billiard map fixes the boundary $\bo\Psi=\bo\Om\times\{0,\pi\}$ pointwise.
We will work with the interior of the phase space, 
and use the notation $\Psi$ for it.

In order to simplify the exposition, in sections~\ref{mirror},~\ref{per},~\ref{calcul},
which form  the basis of the paper, we restrict the setting to  
convex billiard tables. In the proof of our main results in
section~\ref{main}, we extend  the discussion to nonconvex billiard tables. 

Let $T\Psi$ be the tangent bundle. 
The coordinate vector fields $\bo/\bo s,\bo/\bo\al$ yield an explicit isomorphism
$T\Psi=\RR\times\Psi$. The canonical invariant measure satisfies $d\mu=\sin\al\,dsd\al$.
Let $(s_1,\al_1)=F(s,\al)$, and let $DF|_{(s,\al)}:T_{(s,\al)}\Psi\to T_{(s_1,\al_1)}\Psi$ 
be the differential of the billiard ball map. The partial derivatives $\bo s_1/\bo s,
\bo s_1/\bo\al$, etc are straightforward to compute. See \cite{GK}, equation 1.1. The invariance of measure 
$\mu$ corresponds to
\begin{equation}   \label{det_eq}
\det\left(DF|_{(s,\al)})\right)=\frac{\sin\al_1}{\sin\al}.
\end{equation}

\begin{figure}[htbp]
\begin{center}
\input{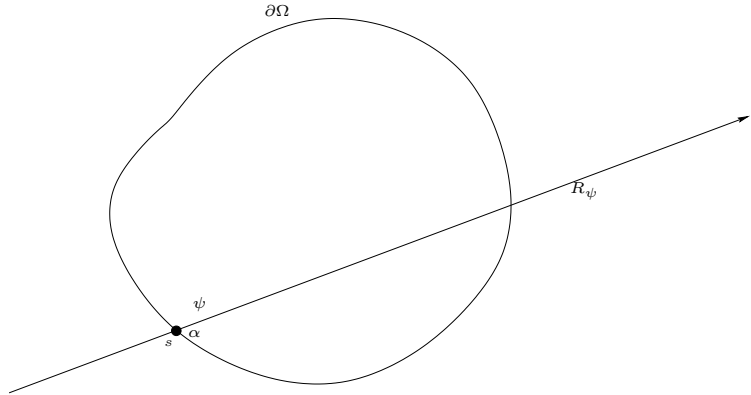}
\caption{The ray $R_{\psi}\subset\RR$ corresponding to a phase point $\psi\in\Psi$.}
\label{ray_fig}
\end{center}
\end{figure}

Let $P(T\Psi)$ be the projectivised tangent bundle. 
Identifying the projective line $P(\RR)$ with the 
augmented real line $\R\cup\infty$, we obtain an isomorphism 
$$
P(T\Psi)=\Psi\times(\R\cup\infty).
$$
We will make this isomorphism explicit by coordinatizing $P(T\Psi)$. 
By a ray we mean an oriented straight line.
For a phase point $\psi\in\Psi$ let $R_{\psi}\subset\RR$ be the corresponding ray. 
See figure~\ref{ray_fig}.
The correspondence $\psi\mapsto R_{\psi}$ yields an isomorphism of
$\Psi(\Om)$ and the space of rays intersecting $\Om$.\footnote{
Here we use the convexity of $\Om$.
}

Next, we augment each $R_{\psi}$ to an oriented projective line by adding the point $\infty$; we will
use the same notation $R_{\psi}$ for the augmented ray. Now we will define a pointwise isomorphism
$P(T_{\psi}\Psi)\sim R_{\psi}$. Vectors in $T_{\psi}\Psi$ correspond to 
infinitesimal deformations of $R_{\psi}$. The rays of the deformation corresponding
to a vector $v\in T_{\psi}\Psi$ intersect $R_{\psi}$ at a point, $f=f(v)\in R_{\psi}$.
The point $f(v)$ does not change under scalings: $f(tv)=f(v),\,t\ne 0$.
On the other hand, any $f\in R_{\psi}$ defines an infinitesimal deformation of $R_{\psi}$,
which is unique, up to scaling. The deformation corresponding to the point $\infty$
moves $R_{\psi}$ parallel to itself. This establishes an isomorphism between two 
projective lines: $R_{\psi}$ and  $P(T_{\psi}\Psi)$.

It will be useful to view elements of $P(T\Psi)$ as pairs $(R_{\psi},f\in R_{\psi})$.
Let $\psi=(s,\al)$, and let $x=x(f)$ be the signed distance from $s$ to $f$;
this identifies $R_{\psi}$ with $\R\cup\infty$. The realisation of  $P(T\Psi)$
as the space of pairs $(R_{\psi},f\in R_{\psi})$, and the signed distance 
yield global coordinates on $P(T\Psi)$.
In this coordinatization, points of the projectivised tangent bundle
$P(T\Psi)$ are the triples $(s,\al,x)$ where $x\in\R\cup\infty$.

Projectivisation of the tangent bundle turns the differential $D\Psi$ into the projectivised 
differential $PDF:P(T\Psi)\to P(T\Psi)$. Let $(s,\al)\in\Psi$ and let 
$F(s,\al)=(s_1,\al_1)$. In our coordinates on $P(T\Psi)$, the mapping
$PDF$ has the global  form $(PDF)(s,\al,x)=(s_1,\al_1,x_1)$.
Since both $x$ and $x_1$ are projective coordinates, and since
the restriction $PDF|_{(s,\al)}:P(T_{(s,\al)}\Psi)\to P(T_{(s_1,\al_1)}\Psi)$ 
is a projective transformation,
the transformation $x_1=x_1(s,\al;x)$ is fractional-linear.
Its coefficients are determined by the geometry of $\bo\Om$ at the points $s$ and $s_1$.

 
%
%
%
\section{Focusing, mirror equation and the projectivised differential}   \label{mirror}
The isomorphism of $R_{\psi}$ and  $P(T_{\psi}\Psi)$ has an intuitive meaning that
comes from the geometric optics. In order to explain this intuition,
we think of the billiard phase points as light rays that are trapped inside the room $\Om$
with a perfectly reflecting wall, $\bo\Om$. Infinitesimal deformations of a light ray $R_{\psi}$
correspond to infinitesimal light beams centered about $R_{\psi}$.
An infinitesimal light beam is determined by its position and by its size.
In the correspondence between light beams and vectors of the tangent space $T_{\psi}\Psi$,
the size of the beam  corresponds to the norm of the tangent vector;
the position of the beam  corresponds to the line in $T_{\psi}\Psi$ spanned by the tangent vector.
Every infinitesimal light beam has a focusing point, $f\in R_{\psi}\cup\infty$,
which determines its position. This is an isomorphism  of 
the projective lines $R_{\psi}$ and  $P(T_{\psi}\Psi)$.

Let $R\in\Psi$ be a light ray, and let $R_1\in\Psi$ be the reflected light ray.
An infinitesimal light beam centered about $R$ and focused at $f\in R$ becomes, after
reflection in $\bo\Om$, an infinitesimal light beam centered about $R_1$.
Let $f_1\in R_1$ be the focusing point of the reflected beam. The transformation
$f\mapsto f_1$ is given by the projectivised differential discussed in section~\ref{set}.
An equation relating $f\in R$ and $f_1\in R_1$ is the {\em mirror equation}
of geometric optics. Let $\ka=\ka(s)$ be the curvature of $\bo\Om$ at the reflection point,
and let $\al$ be the angle of incidence.

\begin{figure}[htbp]
\begin{center}
\input{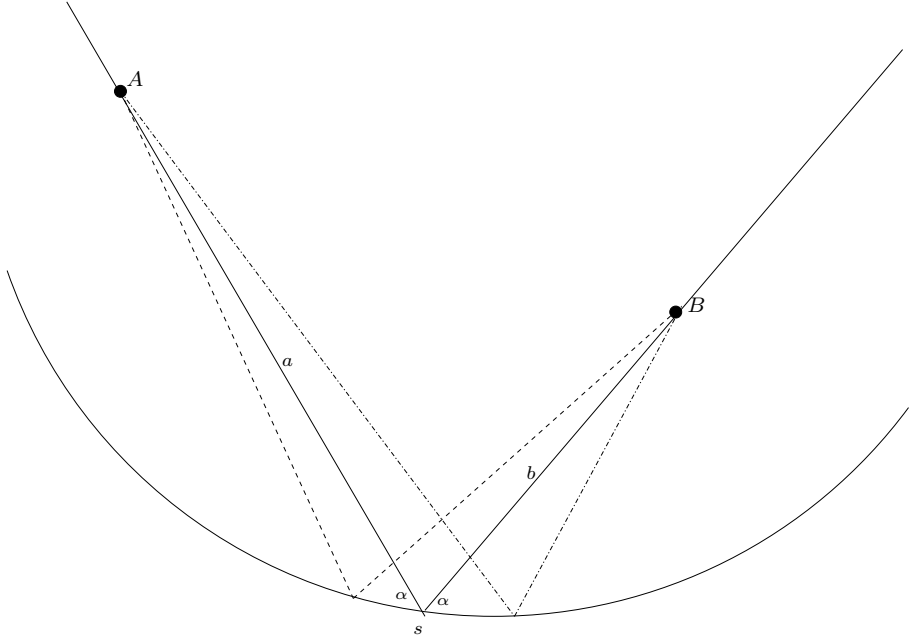}
\caption{Focusing of reflected infinitesimal beams: The mirror equation.}
\label{mirro_eq_fig}
\end{center}
\end{figure}

Let $a$ and $b$ be the signed distances from the focusing points to the reflection point.
See figure~\ref{mirro_eq_fig} for notation. Then the mirror equation says
\begin{equation}  \label{mirr_eq}
\frac{1}{a} + \frac{1}{b}  = \frac{2\ka}{\sin\al}.
\end{equation}
It follows from elementary euclidean geometry. 
See \cite{GK}, Proposition 1.1.

\begin{figure}[htbp]
\begin{center}
\input{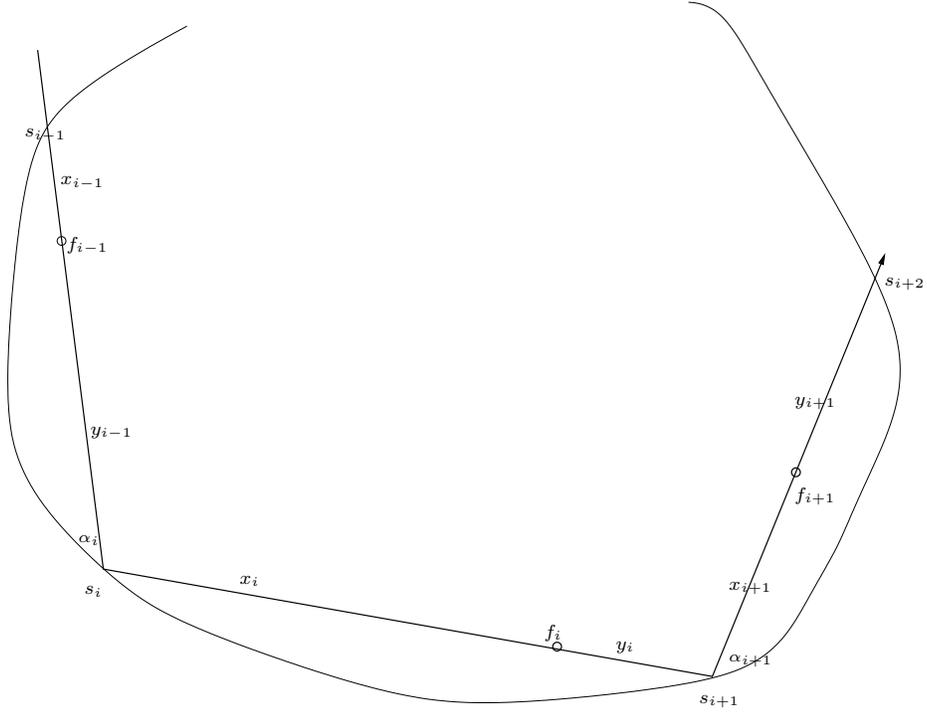}
\caption{An orbit of the projectivised differential of the billiard map.}
\label{focus_fig}
\end{center}
\end{figure}

\medskip
We will view $\Psi\times(\R\cup\infty)$ as an extended phase space, and will
study the orbits of $PDF:\Psi\times(\R\cup\infty)\to\Psi\times(\R\cup\infty)$.
Let $(\psi,x)\in\Psi\times(\R\cup\infty)$ be an extended phase point. Let
$\{(\psi_i,x_i)=(PDF)^i(\psi,x),\,i\in\Z\}$, be the orbit.
See figure~\ref{focus_fig} for notation.

More precisely,
the footpoints are $s_i$; the rays are $R_i$; the lengths are $l_i$;
the $i$th focusing point is $f_i\in R_i$; its coordinate in $R_i$ is $x_i$;
we set $y_i=l_i-x_i$. The lemma below gives explicit expressions for
the fractional-linear transformations $B_{i+1}:P(T_{\psi_i}\Psi)\to P(T_{\psi_{i+1}}\Psi)$
and their inverses $B_i^{-1}:P(T_{\psi_i}\Psi)\to P(T_{\psi_{i-1}}\Psi)$. Set
\begin{equation}  \label{mirr_recur_eq1}
\frac{2\ka_i}{\sin\al_i}=\frac{1}{\la_i}.
\end{equation}
\begin{lem}     \label{mirr_recur_lem}
Let $\bo\Om$ be a curve of class $C^2$ with strictly positive curvature.
Let $(\psi,x)\in\Psi\times(\R\cup\infty)$ be arbitrary. Then for $i\in\Z$ we have
\begin{equation}  \label{mirr_recur_eq2+}
x_{i+1}  = \la_{i+1} - \frac{\la_{i+1}^2}{x_i-(l_i-\la_{i+1})}
\end{equation}
and 
\begin{equation}  \label{mirr_recur_eq2-}
y_{i-1}  = \la_i - \frac{\la_i^2}{y_i-(l_i-\la_i)}.
\end{equation}
\begin{proof}
By equation~\eqref{mirr_eq}, we have (for any $i$)
\begin{equation}  \label{mirr_recur_eq3}
\  \frac{1}{y_i} + \frac{1}{x_{i+1}} =  \frac{1}{\la_{i+1}}; 
\ \ \frac{1}{y_{i-1}} + \frac{1}{x_i} =  \frac{1}{\la_i}.
\end{equation}
Solving for $x_{i+1}$ and using $x_i+y_i=l_i$,
we obtain the recurrence relation equation~\eqref{mirr_recur_eq2+}.
The derivation of equation~\eqref{mirr_recur_eq2-} goes the same way.
\end{proof}
\end{lem}

We will establish conventions relating matrices
$A=(a,b;c,d)\in\gl$ and fractional-linear transformations by $[A]\in\pgl=\psl$.
Often, we will not notationally distinguish between 
fractional-linear transformations and their
representing matrices; if  $A=(a,b;c,d)\in\gl$, we will use the notation
%
$$
A\circ z = \frac{az+b}{cz+d},
$$
%
and call rational functions of this kind {\em fractional-linear}.
Note that a fractional-linear function is linear iff $c=0$.
Thus $c\ne 0$ or $c=0$ is a property of $[A]$.

Let the setting be as in Lemma~\ref{mirr_recur_lem}.
The fractional-linear transformation $x_0\mapsto x_i$ is the {\em focusing on $i$th step}.
Set $A_i=B_i\cdots B_1$. Then the focusing on $i$th step is given by
$$
x_i=A_i\circ x_0.
$$

%
\begin{defin}  \label{non_lin_foc_def}
Let $\psi_0\in\Psi$, and let $A_i=(a_i,b_i;c_i,d_i),\,i\in\Z,$ be as above.
If $c_i=0$ we say that the {\em focusing on $i$th step is linear}.
If $c_i\ne 0$ we say that the {\em focusing on $i$th step is nonlinear}.
\end{defin} 
%


%
\begin{lem}     \label{non_lin_foc_lem}
The focusing on the first step is always nonlinear.
The focusing cannot be linear twice in a row.
\begin{proof}
Since the focusing on $0$th step is linear, the former claim
is a particular case of the latter, which is immediate from equation~\eqref{mirr_recur_eq2+}.
\end{proof}
\end{lem}

\medskip


We will not need explicit formulas for the focusing transformations
$x_i=A_i\circ x_0$. We will use only their derivatives.

\medskip

\begin{lem}     \label{deriv_foc_lem}
Let $\psi_0\in\Psi$, and let $\psi_i=F^i(\psi_0),\,0\le i,$ be the orbit.
Let $x_0\in\pr$ be arbitrary, and let $x_i\in\pr$ be the focusing distance on $i$th step.
Then for $1\le k$ we have
\begin{equation}  \label{der_foc_eq}
\frac{dx_k}{dx_0}= \prod_{i=1}^k \frac{(x_i-\la_i)^2}{\la_i^2}.
\end{equation}
\begin{proof}
Differentiating equation~\eqref{mirr_recur_eq2+}, we obtain
$$
dx_{i+1}=\frac{x_{i+1}^2}{y_i^2}\ dx_i.
$$
Solving equation~\eqref{mirr_recur_eq3} for $y_i$, substituting
the expression into the equation above, and iterating  back to $i=0$, we
obtain the claim.  
\end{proof}
\end{lem}

\section{Periodic points and focusing}  \label{per}
\begin{figure}[htbp]
\begin{center}
\input{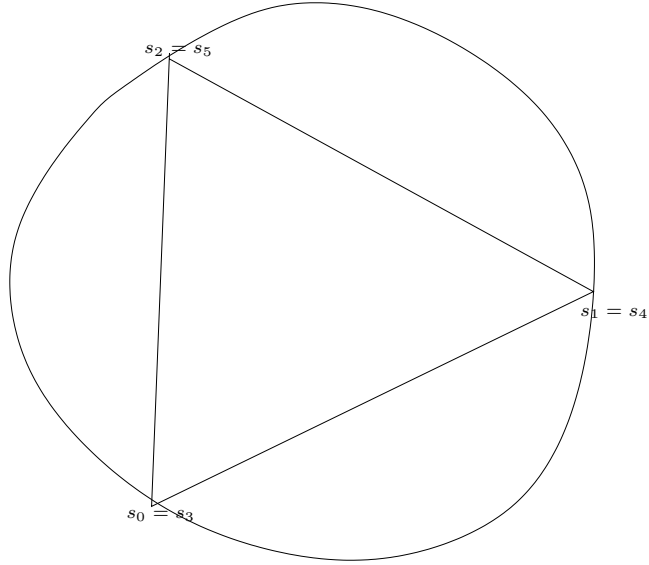}
\caption{A $6$-periodic phase point; the corresponding harmonic polygon
is the triangle $s_0s_1s_2$ traced twice.}
\label{harmo_fig}
\end{center}
\end{figure}

Let $n\in\N$. A phase point $\psi\in\Psi$ is {\em $n$-periodic} if $F^n(\psi)=\psi$.
The smallest such $n$ is the {\em minimal period} of $\psi$.
A $n$-periodic point determines a harmonic $n$-gon,
$P=P(\psi,n)\subset\Om$. We label the quantities pertaining to $P$ by $i\in\Z$ 
with the cyclic convention $i+n=i$, and set $\psi_0=\psi$. If $n$ is not the minimal
period, then $P$ is a polygon with a multiplicity.
See figure~\ref{harmo_fig} for an example.

The differential $(DF_{\psi})^n$ determines the type of periodic point.
By equation~\eqref{det_eq}, $(DF_{\psi})^n\in\mac$. The $n$-periodic point $\psi\in\Psi$ 
is {\em hyperbolic, elliptic, or parabolic} if $(DF_{\psi})^n$ is 
hyperbolic, elliptic, or parabolic respectively.\footnote{
Note that the type of a $n$-periodic point depends both on the point and the period 
which need not be minimal.
}
By our convention, the matrices $\pm1\in\mac$ are not parabolic; if
$(DF_{\psi})^n=\pm1$, then $\psi$ is a {\em degenerate $n$-periodic point}.

\vspace{2mm}

Instead of working with differentials, we will use  the projectivised differentials,
or, equivalently, the focusing points.
Let $f_0\in R_0$ be an arbitrary focusing point,
and let $x_0\in\R\cup\infty$ be the signed focusing distance.
Then $f_0$ generates an infinite sequence, $f_i\in R_i,1\le i,$ of focusing points, and 
let $x_i\in\R\cup\infty$ be the focusing distances. 
For $0\le i$ let $B_i=P(DF_{\psi_i})\in\psl$
be the projectivised differential at $\psi_i$. For $0\le i<j$ set $A_{j,i}=B_{j-1}\cdots B_i$.
%
%
\begin{lem}     \label{neutr_per_prop}
Let $\psi\in\Psi$ be a $n$-periodic point. Then the following
statements are equivalent:

\noindent 1. The point $\psi$ is degenerate;

\noindent 2. The identity $x_n\equiv x_0$ holds;

\noindent 3. For $0\le i \le n-1$ we have $dx_{i+n}/dx_i\equiv 1$.
\begin{proof}
Set $D_i=A_{n+i,i}\in\psl$. Then $\psi$ is degenerate iff $D_0=1$ iff $D_i=1$ for all $i$.
Thus, properties 1 and 2 are equivalent, and 3 follows from 1.
It remains to show that 3 implies 1.

Assume the opposite, i. e., that $\psi$ satisfies condition 3, but it is not degenerate.
Then for $0\le i \le n-1$ we have $x_{i+n} \equiv x_i +b_i$ for some $b_i\ne 0$. 
Equivalently, every $D_i=B_{n+i-1}\cdots B_i$ is unipotent and upper-triangular.
By periodicity, $B_{n+i}=B_i$. Therefore
$$
D_1=B_nB_{n-1}\cdots B_1=B_0\left(B_{n-1}\cdots B_1B_0\right)B_0^{-1}=B_0D_0B_0^{-1}.
$$
Thus, $D_0,D_1$ are unipotent, upper-triangular, and $B_0$ conjugates them.
Therefore $B_0$ is upper-triangular. But, by equation~\eqref{mirr_recur_eq2+}, 
$B_0$ is not upper-triangular.
\end{proof}
\end{lem}
\begin{prop}     \label{neutr_per_cor}
Let $\psi\in\Psi$ be $n$-periodic. Designate
any point of the orbit $\psi,F(\psi),\dots,F^{n-1}(\psi)$ as $\psi_0$;
let $\la_0,\dots,\la_{n-1}$ be the corresponding numbers. 
Let $x_0\in\R\cup\infty$, and let $x_i=x_i(x_0),1\le i \le n-1,$
be the corresponding focusing distances. 
Then $\psi$ is degenerate iff for any choice of $\psi_0$ and arbitrary $x_0$ we have 
\begin{equation}  \label{neutr_per_eq}
 \prod_{i=0}^{n-1}(x_i-\la_i) = \pm\prod_{i=0}^{n-1}\la_i.
\end{equation}
\begin{proof}
By Lemma~\ref{deriv_foc_lem}, equation~\eqref{neutr_per_eq}
is equivalent to condition 3 of Lemma~\ref{neutr_per_prop}.
\end{proof}
\end{prop}

\medskip

Let $\psi\in\Psi$. For $0\le i$ set $\psi_i=F^i(\psi)$.
Let $s_0,s_1,\ldots\in\bo\Om$ be the corresponding sequence of footpoints.
Let $\la_i$ (resp. $l_i=l(s_i,s_{i+1})$) be the corresponding numbers 
(resp. side lengths).

\begin{prop}     \label{neutr_geom_prop}
Suppose that $\psi$ is $n$-periodic.
Then $\psi$ is degenerate iff $n=2m$ is even and the following statements hold:\\

\noindent 1. For all $i$ we have
\begin{equation}  \label{basic_geom_eq}
l_i=\la_i+\la_{i+1};
\end{equation}

\noindent 2. We have
\begin{equation}  \label{prod_eq}
\la_0\la_2\cdots\la_{2m-2}=\la_1\la_3\cdots\la_{2m-1}.
\end{equation}
\begin{proof}
Let $B_i,A_{j,i}$ be as in Lemma~\ref{neutr_per_prop}, and set $A_k=A_{k,0}$.
Then for $0\le i$ the rational function $x_i=A_i\circ x_0$ has one of the expressions
\begin{equation}  \label{frac_lin_eq}
x_i= p_i + q_i(x_0-r_i)^{-1};\ \ x_i= q_ix_0-p_i
\end{equation}
where $p_i,q_i,r_i\in\R$ with $q_i\ne 0$. Thus, $x_i$ has 
the former (resp. latter) form iff the focusing on $i$th step is nonlinear (resp. linear).

The rational function $F(x_0)=\prod_{i=0}^{n-1}(x_i-\la_i)$
is a product of $n$ fractional-linear functions $x_i(x_0)-\la_i$.
Suppose first that $\psi$ is degenerate.
Then, by Lemma~\ref{neutr_per_prop}, $F=const$. 

Let $I\subset\{0,1,\dots,n-1\}$ (resp. $J\subset\{0,1,\dots,n-1\}$) be the set of indices
such that $x_i(x_0)-\la_i$ is linear (resp. nonlinear).
Then the asymptotics of $F(x_0)$ as $x_0\to\infty$ satisfy
$$
F(x_0)=\con\  x_0^{|I|-|J|} + O(x_0^{|I|-|J|-1})
$$
with $\con\ne 0$. Since, by Lemma~\ref{non_lin_foc_lem}, $|I|\le n/2$,
our assumption implies that $|I|=|J|=n/2$, hence $n=2m$. Since $0\in I$,
by Lemma~\ref{non_lin_foc_lem},
$I=\{0,2,\dots,2m-2\},J=\{1,3,\dots,2m-1\}$. Thus, the focusing is linear
on every even step and nonlinear on every odd step.

By equation~\eqref{mirr_recur_eq2+}
\begin{equation}  \label{foc_rec_eq}
x_{i+1}-\la_{i+1}=-\frac{\la_{i+1}^2}{(x_i-\la_i)-(l_i-\la_i-\la_{i+1})}.
\end{equation}
Calculating from equation~\eqref{foc_rec_eq} the expression $x_{i+2}-\la_{i+2}$
in terms of $x_i-\la_i$, we see that the function $x_2(x_0)$ is linear iff
$l_1-\la_1-\la_2=0$.

By periodicity, we can do the preceding analysis starting from any point of the orbit.
Starting from $i-1$, we obtain $l_i-\la_i-\la_{i+1}=0$.
Thus the degeneracy implies the identities equation~\eqref{basic_geom_eq}.

Conversely, suppose that the identities equation~\eqref{basic_geom_eq} hold.
Then by elementary calculations, for $0\le k\le m-1$, we have the focusing equations
$$
x_{2k+1}-\la_{2k+1} = -\left(
\frac{\la_{2k+1}^2\la_{2k-1}^2\cdots\la_1^2}{\la_{2k}^2\la_{2k-2}^2\cdots\la_2^2}
\right)(x_0-\la_0)^{-1};
$$
$$
x_{2k}-\la_{2k} = \left(
\frac{\la_{2k}^2\la_{2k-2}^2\cdots\la_2^2}{\la_{2k-1}^2\la_{2k-3}^2\cdots\la_1^2}
\right)(x_0-\la_0).
$$
Multiplying out these equations, we obtain
$$
\prod_{i=0}^{2m-1}(x_i-\la_i)=(-1)^m\la_1^2\la_3^2\cdots\la_{2m-1}^2.
$$
Compairing this with Lemma~\ref{neutr_per_prop}, we conclude that
equation~\eqref{basic_geom_eq} implies the degenerate focusing iff the identity
equation~\eqref{prod_eq} holds.
\end{proof}
\end{prop}

\begin{rem}   \label{sign_rem}
{\em The argument above shows that the sign in equation~\eqref{neutr_per_eq}
is determined by the pairity of $n/2$. Namely,
case i) holds iff $n=0\mod{4}$; case ii) holds iff $n=2\mod{4}$. 
}
\end{rem}

\section{Degenerate periodic points and a functional  relation}   \label{calcul}
\begin{figure}[htbp]
\begin{center}
\input{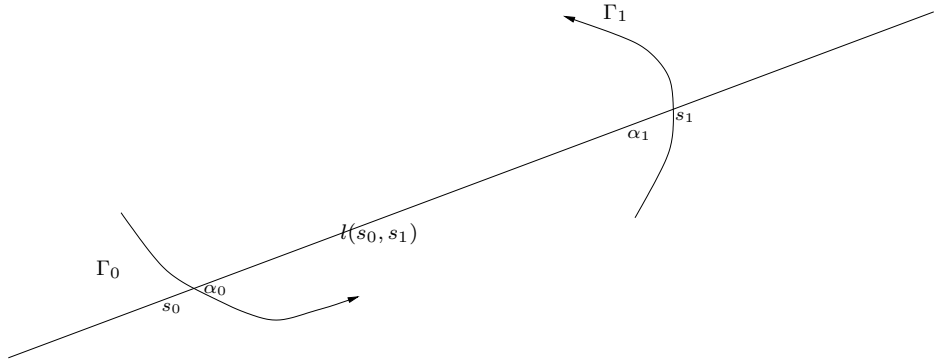}
\caption{A billiard chord and its attributes.}
\label{chord_fig}
\end{center}
\end{figure}

In this section we will study the relation implied by equation~\eqref{basic_geom_eq}. 
We will view phase points $\psi\in\Psi$ as directed chords; they correspond to pairs 
$(s_0,s_1)$. See figure~\ref{chord_fig}.

Let $\ka(s),\rho(s)$ denote the curvature and the radius of curvature.
Let $l(s_0,s_1)$ be the chord length. Differentiating with respect to $s_0$ and $s_1$,
and using elementary geometry, we obtain:
\begin{equation}  \label{chord_der_eq}
dl(s_0,s_1)=-\cos\al_0 ds_0 + \cos\al_1 ds_1,
\end{equation}
and
\begin{equation}  \label{angl_der_eq1}
d\al_0=\left(\frac{\sin\al_0}{l}-\ka(s_0)\right)d\,s_0\ +\ \left(\frac{\sin\al_1}{l}\right)d\,s_1,
\end{equation}
and
\begin{equation}  \label{angl_der_eq2}
d\al_1=
\left(-\frac{\sin\al_0}{l}\right)d\,s_0\ +\ \left(-\frac{\sin\al_1}{l}+\ka(s_1)\right)d\,s_1.
\end{equation}

In the proposition below by a curve we mean a smooth immersion
$f:I\to\RR$, where $I\subset\R$ is an open interval. We identify
a curve with its image $\Ga\subset\RR$. By the smoothness class
of $\Ga$ we will mean the smoothness  of the mapping $f:I\to\Ga$.

\begin{prop}   \label{gener_anal_prop}
Let $\Ga_0,\Ga_1\subset\RR$ be $C^3$ curves and let $c$ be a constant. Suppose  that 
for all $s_i\in\Ga_i$ we have
\begin{equation}  \label{gener_anal_eq}
\rho(s_0)\sin\al_0 + \rho(s_1)\sin\al_1 = (1+c)l(s_0,s_1).
\end{equation}
Then  $c=0$ and both $\Ga_0,\Ga_1$ are arcs of a circle.
\begin{proof}
It is classical that the identity equation~\eqref{gener_anal_eq} with $c=0$
holds for any pair of arcs of a circle.


Differentiating the identity equation~\eqref{gener_anal_eq},
and using relations~\eqref{chord_der_eq},~\eqref{angl_der_eq1},~\eqref{angl_der_eq2},
we obtain
\begin{equation}  \label{ident_eq0}
\rho'(s_0)+\frac{\cos\al_0}{l}\rho(s_0)-\frac{\cos\al_1}{l}\rho(s_1)=-c\cot\al_0, 
\end{equation}
and
\begin{equation}  \label{ident_eq1}
\rho'(s_1)-\frac{\cos\al_1}{l}\rho(s_1)+\frac{\cos\al_0}{l}\rho(s_0)=c\cot\al_1. 
\end{equation}
Differentiating the identity equation~\eqref{ident_eq0} with respect to $s_1$,
using the relations equations~\eqref{angl_der_eq1},~\eqref{angl_der_eq2},
and solving for $\rho'(s_1)$, we obtain the expression
\begin{equation}  \label{ident_eq2}
\rho'(s_1)=-\frac{\cos(\al_0-\al_1)}{l\cos\al_1}\rho(s_0)+\frac{\cos 2\al_1}{l\cos\al_1}\rho(s_1)
+\tan\al_1\left[1-\frac{c}{\sin^2\al_0}\right].
\end{equation}

Subtracting equation~\eqref{ident_eq2} from equation~\eqref{ident_eq1}, we obtain
a linear relation between $\rho(s_0),\rho(s_1)$. After simplifications, we have
$$
\frac{\sin\al_0\tan\al_1}{l}\rho(s_0)+\frac{\sin\al_1\tan\al_1}{l}\rho(s_1)=
\tan\al_1-c\left[\cot\al_1+\frac{\tan\al_1}{\sin^2\al_0}\right].
$$
A further simplification yields
\begin{equation}  \label{linear_eq}
\tan\al_1\left[\frac{\sin\al_0\rho(s_0)+\sin\al_1\rho(s_1)}{l}\right]=
\tan\al_1-c\left[\cot\al_1+\frac{\tan\al_1}{\sin^2\al_0}\right].
\end{equation}
From equations~\eqref{gener_anal_eq},~\eqref{linear_eq}, we have
$$
1+c=1-c\left[\cot^2\al_1+\frac{1}{\sin^2\al_0}\right]
$$
implying $c=0$.

\medskip

It remains to show that the identity $\sin\al_0\rho(s_0)+\sin\al_1\rho(s_1)=l(s_0,s_1)$
implies that $\Ga_0,\Ga_1$ are arcs of the same circle.
Subtracting equation~\eqref{ident_eq1} from equation~\eqref{ident_eq0}, we obtain
the identity
$$
\rho'(s_1)-\rho'(s_0)=c\left[\cot\al_0+\cot\al_1\right].
$$
Since $c=0$, by our preceding calculations, this implies $\rho=\const$.
\end{proof}
\end{prop}

Note that Proposition~\ref{gener_anal_prop} yields a characterization of euclidean
circles. Namely, a round circle is singled out by an identity
satisfied by the distance between points of any pair of its arcs, no matter how
small they are.
The following statement is immediate from Proposition~\ref{gener_anal_prop}.

\medskip

\begin{corol}   \label{no_iden_cor}
Let $\Ga_0,\Ga_1\subset\RR$ be $C^3$ curves. For $s_i\in\Ga_i$
let $\al_i$ be the respective angles, and let $l(s_0,s_1)$ be the
distance. The relation  
\begin{equation}  \label{imposs_eq}
\rho(s_0)\sin\al_0 + \rho(s_1)\sin\al_1 = 2\,l(s_0,s_1)
\end{equation}
cannot hold identically.
\end{corol}

\section{Main results}   \label{main}
We will now apply the preceding material to study periodic points for the billiard
in a planar, compact domain, with a piecewise $C^3$ boundary.
Denote by $\el(\Om),\hy(\Om),\pa(\Om),\neut(\Om)\subset\per(\Om)\subset\Psi$
the sets of elliptic, hyperbolic, parabolic and degenerate periodic points respectively.
The subscript indicates the period; for instance, $\el_n(\Om)$ is the
set of elliptic $n$-periodic points.

%
%
%
%
\begin{thm}  \label{emt_inter_thm}
Let $\Om\subset\RR$ be a compact domain with a piecewise $C^3$ boundary.
Then the set of periodic points for the billiard in $\Om$ has an empty interior.
\begin{proof}
Assume first that $\Om$ is convex, and that $\bo\Om$ is $C^3$.
Suppose that the claim is false, and let $O\subset\Psi$ be an open set of
periodic points. Let $F:\Psi\to\Psi$ be the billiard map.
By taking a smaller $O$, if need be, we obtain $n\in\N$ such that
$F|_O^n=\id$. Thus, all $\psi\in O$ are degenerate $n$-periodic points.

Let $s_0(\psi),s_1(\psi)$ be the endpoints of the chord of $\Om$ defined by $\psi$.
As $\psi$ varies in $O$, the points $s_0(\psi),s_1(\psi)$ run through 
open arcs $\Ga_0,\Ga_1\subset\bo\Om$. By equation~\eqref{basic_geom_eq},
the pairs $s_0,s_1$ satisfy the identity equation~\eqref{imposs_eq}.
By Corollary~\ref{no_iden_cor}, this is impossible. This establishes the claim for
convex $\Om$. 

We observe that convexity of $\Om$ has not been used in the analysis of the
projectivised differential that led to equation~\eqref{basic_geom_eq}.
It was convenient to assume convexity in order to identify phase points
with rays intersecting $\Om$. A local version of this identification
does not require convexity. In some of the analytic arguments that
yielded equation~\eqref{basic_geom_eq} we used the positivity of $\la=\rho(s)\sin\al/2$.
These arguments remain valid as long as $\rho(s)$ does not change sign.
The change of sign can occur only at points of $\bo\Om$ where its curvature vanishes, becomes
infinite, or jumps. The two latter possibilities can happen only at isolated points
of $\bo\Om$. By going to a smaller open set in $\Psi$, if need be, we can ignore these
points. The same argument holds if the curvature vanishes at isolated points.
Vanishing of the curvature on an interval means that  $\bo\Om$ contains
a straight arc. If the set of footpoints of $\psi\in O$ contains any curved arcs,
then we get rid of flat arcs via the reflection trick \cite{Gut03}.
This trick changes $\bo\Om$; the curve that replaces $\bo\Om$ has an open set of
periodic billiard orbits which do not pass through flat arcs. 
Thus, we can apply the preceding argument. If the billiard orbits generated by $O$
encounter only flat arcs in  $\bo\Om$, then we have a polygonal billiard table
with an open set of periodic orbits. This is impossible \cite{Gut03}.

We summarise the preceding argument as follows. Assuming an open set of
periodic billiard orbits, we either come to a contradiction right away,
or arrive to  a situation considered in Corollary~\ref{no_iden_cor}. 
But this is a contradiction as well. 
\end{proof}
\end{thm}

Several publications have proved the nonexistence of an open set of $3$-periodic
billiard orbits \cite{BZh,R89,St,Vo97,W94}. Recall that periodic billiard points 
that belong to an open set of those are necessarily degenerate. 
The converse fails: A billiard table may have isolated
degenerate periodic orbits.

\begin{thm}  \label{odd_neut_type_thm}
Let $\Om\subset\RR$ be a compact domain with a piecewise $C^3$ boundary.
Then for any odd $n$ there are no degenerate $n$-periodic points in $\Om$.
\begin{proof}
For convex $\Om$ this is a part of Proposition~\ref{neutr_geom_prop}.
The argument in the proof of Theorem~\ref{emt_inter_thm} 
that shows how to remove the assumption 
of convexity in that statement, works here just the same. 
\end{proof}
\end{thm}

A lebesgue-type measure on a manifold is a measure that is given by
a local differentiable density. The canonical invariant measure for the billiard
in $\Om\subset\RR$ is lebesgue-type.

\begin{thm}   \label{main_conv_thm}
Let $\Om\subset\RR$ be a compact domain with a piecewise $C^3$ boundary.
Let $\mu$ be a lebesgue-type measure on $\Psi$. Then $\mu(\per(\Om))=0$.
\begin{proof}
Let $n\in\N$ be arbitrary. By Theorem~\ref{per_crit_thm} of section~\ref{crit},  
the set $\per_n(\Om)\subset\Psi$ of 
$n$-periodic points is a finite union of submanifolds of smaller dimension. 
\end{proof}
\end{thm}

\noindent

\medskip

%


\section{Appendix 1: Periodic points of the billiard map and critical locus
of the perimeter function}  \label{crit}
The set of $n$-periodic points of a conservative surface diffeomorphism,
$F:S\to S$, may have a complicated structure. In particular, it may be nowhere dense
and have positive measure.
We will outline an example.

Let $T^2=S^1\times S^1$ be the standard torus, and let $d\mu=dx\,dy$ be the lebesgue measure.
We consider selfmappings of the form $F(x,y)=(x+f(y),y),\,G(x,y)=(x,y+g(x))$.
Then $F:T^2\to T^2$, $G:T^2\to T^2$) are 
conservative diffeomorphisms of the same smoothness as the functions 
$f:S^1\to S^1$, $g:S^1\to S^1$). The composition $H=G\circ F:T^2\to T^2$
satisfies  $H(x,y)=(x+f(y),y+g(x+f(y))$. Let $Z_f,Z_g\subset S^1$ be the zero sets of $f,g$, 
and let $\fix(\cdot)$ denote the set of fixed points of a diffeomorphism. 
Then $\fix(H)=Z_f\times Z_g$.

For any kantor set $K\subset S^1$ there exists a $C^{\infty}$  function $f:S^1\to S^1$
such that $K=Z_f$. Combining this remark with the preceding analysis, we obtain the following
statement. Let $K',K''\subset S^1$ be arbitrary cantor sets.
Then there exists a conservative, $C^{\infty}$ diffeomorphim $H:T^2\to T^2$ such that
$\fix(H)=K'\times K''$.

This yields examples of conservative surface diffeomorphisms with 
fractal sets of periodic points.
However, this cannot happen for the billiard, in view of a relation between
its set of $n$-periodic points 
and the critical locus of the perimeter function 
on  $(\bom)^n$. We will now establish this relation.

\medskip 

Let $\Om\subset\RR$ be a compact, strictly convex\footnote{
We make this assumption for simplicity of exposition.
}
domain with a $C^k$ boundary. Rescaling the metric if need be, we assume without
loss of generality that $|\bom|=1$, and identify $\bom^n$ with the torus
$T^n=\times^n(S^1)=\{(s_0,\ldots,s_{n-1})\}$. Set
\begin{equation}  \label{perim_eq}
p(s_0,\ldots,s_{n-1})=l(s_0,s_1)+\cdots+l(s_{n-1},s_0).
\end{equation} 
Note that $p(s_0,\ldots,s_{n-1})$ is the perimeter of the inscribed ``polygon"
$P(s_0,\ldots,s_{n-1})$ with corners $s_0,\ldots,s_{n-1}$. (By
parentheses, we indicate that consecutive corners $s_i,s_{i+1}$ need not be distinct.
When this is the case, we will say that $P(s_0,\ldots,s_{n-1})$ is a degenerate polygon.
Otherwise, $P(s_0,\ldots,s_{n-1})$ is a nondegenerate polygon.)

\medskip

By a $C^k$ curve in $\Psi$ we will mean the range of a $C^k$ immersion $f:[0,1]\to\Psi$.

\begin{thm}  \label{per_crit_thm} 
Let $\Om\subset\RR$ be a compact domain with a piecewise $C^3$ boundary.
Then for any $n\in\N$ the set $\per_n(\Om)\subset\Psi$ is a finite union of $C^2$
curves and points. 
\begin{proof}
Denote by $q_n:T^n\to\Psi$ the projection that sends $(s_0,\ldots,s_{n-1})$
to $(s_0,s_1)\in\Psi$. We will deduce the statement from a sequence of claims.

\noindent 1. Let $\crit_n(\Om)\subset T^n$ be the critical locus 
of the function given by equation~\eqref{perim_eq}. Then for $(s_0,\ldots,s_{n-1})\in\crit_n(\Om)$
the ``polygon" $P(s_0,\ldots,s_{n-1})$ is either a harmonic $n$-gon
or $s_0=\dots=s_{n-1}$. 

A string of indices $0\le i<j \le n-1$ is degenerate
if $s_i=\dots=s_j$.
Let $0\le i<j \le n-1$ be a maximal degenerate string in 
$(s_0,\ldots,s_{n-1})\in T^n$. If $0<i$ or $j<n-1$, then, by equation~\eqref{chord_der_eq},
$(s_0,\ldots,s_{n-1})\notin\crit_n(\Om)$. This proves claim 1.

\medskip

We identify the set of completely  
degenerate $n$-tuples $(s_0,\ldots,s_{n-1})$ with
$\bom$, and denote by $\h_n(\Om)$ the set of harmonic $n$-gons.
Thus $\crit_n(\Om)=\h_n(\Om)\cup\bom$, a disjoint union. 

\noindent 2. The restriction of $q_n:T^n\to\Psi$ to $\h_n(\Om)$ yields a surjective, one-to-one 
mapping $q_n:\h_n(\Om)\to\per_n(\Om)$. Let $(s_0,\ldots,s_{n-1})\in\crit_n(\Om)$.

The criticality condition
equation~\eqref{chord_der_eq} means that $(s_1,s_2)=F(s_0,s_1)$, etc.
Claim 2 is immediate from this observation.

\medskip

\noindent 3. The restriction $q_n:\h_n(\Om)\to\per_n(\Om)$ is a $C^2$ diffeomorphism.

By claim 2 and compactness of $\h_n(\Om)$, the map $q_n:\h_n(\Om)\to\per_n(\Om)$ 
is a homeomorphism. Let $\psi=(s_0,s_1)\in\Psi$ be arbitrary.
The point $(s_0,\ldots,s_{n-1})\in  T^n$ is obtained by applying $F^i$ to $\psi$.
Hence the smoothness class of the mapping $r_n:\Psi\to T^n$ given by 
$r_n(s_0,s_1)=(s_0,\ldots,s_{n-1})$ is that of $F$; the latter is $1$ less than
the smoothness of $\bom$. By our assumption, it is $C^2$. By claim 1, the restriction of
$r_n$ to $\per_n(\Om)$ is the inverse of $q_n$. The claim follows.

\medskip

The set $\h_n(\Om)\subset T^n$ is (essentially) the critical locus of
the function $p_n$ given by equation~\eqref{perim_eq}. 
By a separate argument, its critical locus is a finite union of $C^2$
submanifolds. By claim 3, these are also submanifolds of $\Psi$, hence
their dimension is at most $2$. But, by Theorem~\ref{emt_inter_thm},
their dimension is less than $2$, i. e., they are points and curves.
\end{proof}
\end{thm}

The following proposition slightly extends Theorem~\ref{per_crit_thm}.

\begin{corol}  \label{per_str_cor}
Let $\Om\subset\RR$ be a compact domain with a piecewise $C^3$ boundary.
Let $n\in\N$.

Then $|\el_n(\Om)\cup\hy_n(\Om)|<\infty$;
the set $\pa_n(\Om)\cup\neut_n(\Om)$ is a finite union of $C^2$
curves and points. 
\begin{proof}
Immediate from Theorem~\ref{per_crit_thm} and the observation that
every point in $\el_n(\Om)\cup\hy_n(\Om)$ is isolated.
\end{proof}
\end{corol}

\section{Appendix 2: Spherical billiard tables with open sets of periodic orbits}\label{rem}
\begin{figure}[htbp]
\begin{center}
\input{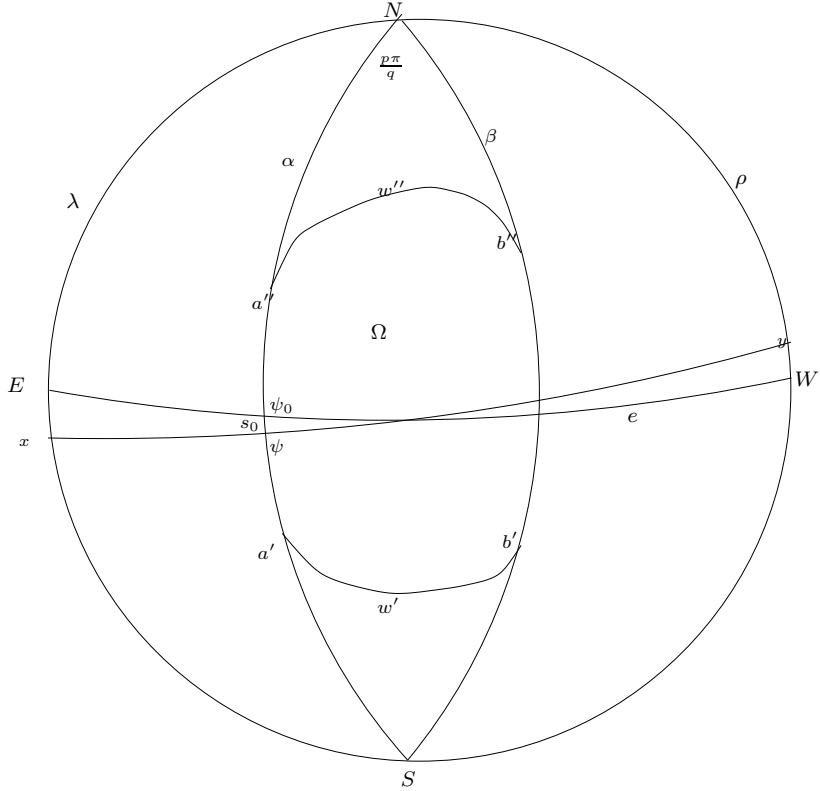}
\caption{A convex billiard table on the round sphere
with an open set of periodic orbits.}
\label{sphere_fig}
\end{center}
\end{figure}

We will outline a construction of compact, $C^{\infty}$ billiard tables $\Om\subset\sph$ 
that have open sets of
periodic orbits. Let $\al,\be\subset\sph$ be two meridians forming a bigon, $B=B(p\pi/q)$, with angle
$p\pi/q$;\footnote{
Here $p,q\in\N$ are relatively prime integers.
} 
let $e$ be the equator. Let $a'\in\al,b'\in\be$ (resp. $a''\in\al,b''\in\be$) be arbitrary points
that are below (resp. above) the equator. Let $w',w''\subset\sph$ be convex arcs 
joining $\al$ and $\be$ tangentially, whose tangency points are $a',b'$ and $a'',b''$ respectively.
Let $\Om\subset\sph$ be the convex domain formed by $\al,\be,w',w''$. 
See figure~\ref{sphere_fig} for notation.

The concepts and the notation that we have developed for the planar billiard pertain to the present
situation as well; we will use them from now on. Let $s_0\in\bom$ be the point where
the equator crosses $\al$. Let $\psi_0=(s_0,\pi/2)\in\Psi$ be the phase point
formed by $e$ at the crossing.

\begin{prop}   \label{sph_ex_prop}
There is an open set $O\subset\Psi$ containing  $\psi_0$ such that $O\setminus\{\psi_0\}$ 
consists of periodic points with rotation number $p/q$.
\begin{proof}
Let $E,W$ be the leftmost, rightmost points of $e$, and let $N,S$ be the north, south poles
respectively. Let $\la,\rho$ be the left, right half-geodesics forming the geodesic in $\sph$
passing through $N,E,S,W$. A geodesic $\ga$ in a neighborhood of $e$ is determined by the points 
$x=x(\ga)=\ga\cap\la,y=y(\ga)=\ga\cap\rho$. Let $\ga(x,y)\subset\sph$ be the corresponding geodesic.

For $0<\ep$ let $G_{\ep}$ be the set of geodesics $\ga(x,y)$ such that $x,y$ are $\ep$-close to $E,W$
respectively.
If $\ep$ is sufficiently small, then any $\ga(x,y)\in G_{\ep}$ intersects $\al$ between $a'$ and $a''$
defining a phase point $\psi(x,y)\in\Psi$. 

We will now compare the billiard orbits in $\Om$ and the bigon $B=B(p\pi/q)$. 
Let $\tPsi$ be the phase space  for the billiard in $B$. 
Using the notation above, let $\tpsi(x,y)\in\tPsi$ be the phase point determined by
the geodesic $\ga(x,y)$ when it intersects $\al$.
We recall the device of {\em unfolding a billiard orbit} defined for any spherical polygon \cite{GT06}. 
In our case it associates with any phase point $\tpsi\in\tPsi$ 
a geodesic, $\ga(\tpsi)$, in $\sph$ and a sequence $B_i$ of isometric images of  $B$ arranged along
$\ga(\tpsi)$, so that any two consecutive images are related by a geodesic reflection. 
We have $\ga(\tpsi(x,y))=\ga(x,y)$.
Using the unfolding, we show that every $\tpsi(x,y)$ is a periodic phase point \cite{GT06};
the sequence of domains $B_i$ returns to the ``initial position" after $2q$ consecutive reflections,
producing a cyclic string of bigons $B_i,\,0\le i \le 2q-1$, with $B_{2q}=B_0$.\footnote{
With the exception of $\tpsi_0$ which is always $2$-periodic \cite{GT06}.
}

Let $\ga_i(x,y)\subset\ga(x,y)$ be the segment of $\ga(x,y)$ inside $B_i$. Note that the unfolding procedures
for the billiard in $\Om$ and the billiard in $B$ coincide as long as  $\ga_i(x,y)$ intersects the
walls  of $B_i$ between $a',a''$ on one side and between $b',b''$ on the other.
The deviation of  $\ga(x,y)$ from the equator is a
continuous function of $x,y$. Therefore for any $a',a'',b',b''$ there exists $0<\ep$
such that every $\ga\in G_{\ep}$ will stay within the requred bounds.

Let $O(\ep)\subset\Psi$ (resp. $\tO(\ep)\subset\tPsi$) 
be the set of phase points corresponding to $\ga(x,y)\in G_{\ep}$.
Since each point $x,y$ varies independently of the other in an open interval, 
we conclude that  $O(\ep)$ (resp. $\tO(\ep)$) is an open neighborhood of $\psi_0$ (resp. $\tpsi_0$).
The unfoldings of orbits of phase points $\psi\in O(\ep)$ and  $\tpsi\in \tO(\ep)$
coincide. The claim now follows from the corresponding claim for $B(p\pi/q)$ \cite{GT06}.
\end{proof}
\end{prop}

For integers $1\le p < q$ denote by $e(p/q)\subset e$ a segment of length $p\pi/q$.
The following assertion is immediate from the preceding proof.
\begin{corol}  \label{sph_ex_cor}
Let $1\le p < q$ be relatively prime. There exist convex $C^{\infty}$ domains $\Om\subset\sph$ 
that contain $e(p/q)$, approximating it arbitrarily closely, and such that the following holds:
i) The billiard map for $\Om$ has an open set $O$ of $q$-periodic points (and $q$ is the minimal period)
with rotation number $p/q$; ii) The corresponding billiard curves have length $2p\pi$.    
\end{corol}

\begin{rem}  \label{sph_ex_rem}
{\em It is worth pointing out that the phenomenon observed in
the examples of Corollary~\ref{sph_ex_cor} holds for nonconvex domains as well.
Indeed, let $\Om$ be as in Corollary~\ref{sph_ex_cor}, and let $\Om\subset\Om'\subset\sph$
be arbitrary domain satisfying
$\bom'\cap\al=\bom\cap\al,\bom'\cap\be=\bom\cap\be$.
Then the billiard in $\Om'$ has an open set of $q$-periodic points. 

Indeed, it is immediate from the preceding discussion that the orbits of phase points
$\psi\in O$ encounter only the common part of $\bom$ and $\bom'$.
Thus, these orbits are ``unaware'' that they ``live" in the bigger billiard table $\Om'$. 
}
\end{rem}

\medskip

%
%

\vspace{3mm}



\begin{thebibliography}{99}
%

\bibitem{BZh} Yu. Baryshnikov and V. Zharnitsky, {\em Sub-Riemannian geometry and 
periodic orbits in classical billiards}, Math. Res. Lett. {\bf 13} (2006), 587 -- 598.




\bibitem{Bi} G.D. Birkhoff, {\em Collected mathematical papers}, Dover Publications, New York 1968. 












\bibitem{Gut01} E. Gutkin, {\em Problems on billiards},\hfill\\ 
http://www.math.iupui.edu/~mmisiure/open/(2001).   
                                                         
\bibitem{Gut03} E. Gutkin, {\em Billiard dynamics: A survey with the emphasis 
on open problems}, Reg. \& Chaot. Dyn. {\bf 8} (2003), 1 -- 13.





\bibitem {GK} E. Gutkin and A. Katok, {\em Caustics for inner and outer billiards}, 
Comm. Math. Phys. {\bf 173} (1995), 101 -- 133.
      

\bibitem {GT06} E. Gutkin and S. Tabachnikov, {\em Complexity of piecewise convex 
transformations in two dimensions, with applications to polygonal billiards 
on surfaces of constant curvature}, Mosc. Math. J. {\bf 6} (2006), 673 -- 701.

\bibitem {Iv80} V. Ivrii, {\em The second term of the spectral asymptotics 
for a Laplace-Beltrami operator on manifolds with boundary}, 
Func. Anal. Appl. {\bf 14} (1980), 98 -- 106. 












%


\bibitem{PeSt} V.M. Petkov and L.N. Stoyanov, {\em Geometry of reflecting rays and inverse
spectral problems}, John Wiley, Chichester 1992.


\bibitem{R89} M.R. Rychlik, 
{\em Periodic points of the billiard ball map in a convex domain}, J. Diff. Geom. {\bf 30} (1989), 
191 -- 205. 





\bibitem{SaVa} Yu. Safarov and D. Vassiliev. {\em The asymptotic distribution of eigenvalues of
partial differential operators}, AMS, Providence 1997.

\bibitem{St} L. Stojanov, 
{\em Note on the periodic points of the billiard}, 
J. Diff. Geom. {\bf 34} (1991), 835 -- 837. 


\bibitem{Tab} S. Tabachnikov, {\em Billiards}, Soc. Math. France, Paris 1995.

 
\bibitem{Vo97} Ya.B. Vorobets, {\em On the measure of the set of periodic points of a billiard}, 
Math. Notes {\bf 55} (1994), 455 -- 460.



\bibitem{WeA} H. Weyl, {\em Gesammelte Abhandlungen}, Springer-Verlag, Berlin 1968.

\bibitem{W94}  M.P. Wojtkowski, 
{\em Two applications of Jacobi fields to the billiard ball problem}, J. Diff. Geom. {\bf 40} (1994), 
155 -- 164. 




\end{thebibliography}
\end{document}